\documentclass{article}
\usepackage{amssymb}
\usepackage{amsmath}
\def\R{{\mathbb{R}}}

\def\Z{{\mathbb{Z}}}
\def\S{{\mathbb{S}}}
\def\w{\mathop{\mathrm{w}}\nolimits}
\def\supp{\mathop{\mathrm{supp}}\nolimits}
\def\wf{\mathop{\mathrm{WF}}\nolimits}
\def\fs{\mathop{\mathrm{FS}}\nolimits}
\def\rh{\mathop{\mathrm{rh}}\nolimits}
\def\sp{\mathop{\mathrm{SP}}\nolimits}
\def\free{\mathop{\mathrm{free}}\nolimits}
\newtheorem{theorem}{Theorem}[section]
\newtheorem{proposition}[theorem]{Proposition}
\newtheorem{lemma}[theorem]{Lemma}
\newtheorem{corollary}[theorem]{Corollary}
\newtheorem{definition}[theorem]{Definition}
\numberwithin{equation}{section}

\title{Singularities of solutions to Schr\"odinger equation
on scattering manifold}
\author{Kenichi {\scshape Ito}\footnote{JSPS Research Fellow,
Graduate School of Mathematical
Sciences, University of Tokyo, 3-8-1 Komaba, Meguroku, Tokyo,
Japan 153-8914: E-mail: \texttt{ito@ms.u-tokyo.ac.jp}}
\ and Shu {\scshape Nakamura}\footnote{%
Graduate School of Mathematical
Sciences, University of Tokyo, 3-8-1 Komaba, Meguroku, Tokyo,
Japan 153-8914: E-mail: \texttt{shu@ms.u-tokyo.ac.jp}. Partially supported by JSPS
Grant Kiban-B 17340033}
}

\date{}

\begin{document}
\maketitle

\begin{abstract}
In this paper we study microlocal singularities of solutions to Schr\"odinger 
equations on scattering manifolds, i.e., noncompact Riemannian manifolds with asymptotically 
conic ends. We characterize the wave front set of the solutions in terms of the initial condition 
and the classical scattering maps under the nontrapping condition. 
Our result is closely related to a recent work by Hassell and Wunsch, 
though our model is more general and the method, which relies heavily on scattering 
theoretical ideas, is simple and quite different.  In particular, we use Egorov-type argument in the 
standard pseudodifferential symbol classes, and avoid using Legendre distributions.
In the proof, we employ a microlocal smoothing property in terms of the \textit{radially 
homogenous wave front set}\/, which is more precise than the preceding results. 
\end{abstract}

\section{Introduction}

We consider Schr\"odinger operators on noncompact manifolds, which is called 
{\em scattering manifolds}\/ following R. Melrose (\cite{Melrose}). 
Scattering manifold is a natural generalization of asymptotically Euclidean manifolds, 
though it is defined as a manifold with boundary with Riemannian metric of a special form. 

Let $M$ be a manifold, and we suppose $M$ is the interior of
a compact manifold $\bar{M}$ with smooth boundary $\partial M$, i.e., 
$M=(\bar{M})^\circ$, where $E^\circ$ denotes the interior of $E$. 
Let $\Omega$ be a neighborhood of $\partial M$ in $\bar{M}$
that is diffeomorphic to $(0,1]\times \partial M$,
where $\{1\}\times \partial M$ corresponds to $\partial M$.
By a map: 
\[
\tau\mapsto r=(1-\tau)^{-1}, \quad (0,1] \stackrel{\cong}\to [1,\infty)
\]
we identify
\begin{align*}
M_\infty:=M\cap \Omega\stackrel{\cong}{\to}(1,\infty)
\times \partial M.
\end{align*}
Such a diffeomorphism is called \textit{a boundary decomposition}.
In this paper we fix a boundary decomposition. 
We also use an atlas of $M_\infty$ such that 
each coordinate neighborhood is diffeomorphic to
$(1,\infty)\times U$, where $U$ is a coordinate
neighborhood of $\partial M$.
We write
\[
M=M_0\cup M_\infty, 
\]
where $M_0\Subset M$ is a relatively compact open submanifold.

\begin{definition}
Let $g^0$, $g$ be Riemannian metrics on $M$.
$g^0$ on $M$ is called \textit{conic} (with respect to a
given boundary decomposition), if there are $R>0$ and
a Riemannian metric $h$ on $\partial M$ such that
\begin{align*}
g^0= dr^2+r^2h_{jk}(\theta)d\theta^jd\theta^k
\quad \mbox{ for } (r,\theta)
\in [R,\infty)\times \partial M.
\end{align*}
A Riemannian metric $g$ is called \textit{scattering metric of long-range type},
if we can write $g=g^0+m$, where
$g^0$ is a conic metric, and
$m$ is of the form
\begin{align*}
m=m^0(r,\theta)dr^2+rm_{j}^1(r,\theta)(drd\theta^j+d\theta^j dr)
+r^2m_{jk}^2(r,\theta)d\theta^jd\theta^k
\end{align*}
such that for some constants $\mu_l>0$ and any indices $k$ and $\alpha$,  
\begin{align*}
|\partial_r^k\partial_\theta^\alpha m^l(r,\theta)|
\le C_{Kk\alpha }r^{-\mu_l-k}
\quad \mbox{ for } (r,\theta)\in (1,\infty)\times K,
\quad l=0,1,2.
\end{align*}
Here $K$ is any compact subset
of a coordinate neighborhood $U\subset \partial M$.
(The dependence of constants $C_{Kk\alpha}$ on $K$ will be suppressed, 
when there is no confusion.)
In particular, if $\mu_0>1$ and $\mu_1>1/2$,
then $g$
is said to be of \textit{short-range type}.
\end{definition}

In this paper we assume $g$ is a scattering metric of 
short-range type.

We consider solutions to the Schr\"odinger equation: 
\begin{align}
i\frac{\partial}{\partial t} u(t,x) = (-\triangle_g+V)u(t,x),\quad
u(0,\cdot)=u_0\in {\mathcal H},
\label{7.8.15.10.0}
\end{align}
where  
\begin{align*}
\triangle_g=\frac{1}{\sqrt{g}}\sum_{j,k=1}^n
\partial_jg^{jk}\sqrt{g}
\partial_k, \quad (g^{jk})=(g_{jk})^{-1},\quad g=\det (g_{jk}),
\end{align*}
is the Laplace-Beltrami operator, and the function space is given by 
\begin{align*}
{\mathcal H}=(L^2(M;\sqrt{g}dx),(\cdot,\cdot)_{{\mathcal H}}),\quad
(u,v)_{{\mathcal H}}=\int_M u(x)\overline{v(x)}\sqrt{g(x)}dx.
\end{align*}
We suppose $V$ is of \textit{smooth short-range type}\/ in
the sense that
$V\in C^\infty(M;\R)$
and for some $\mu_3> 1$
\begin{align*}
|\partial_r^k\partial_\theta^\alpha V(r,\theta)|
\le C_{k\alpha}r^{2-\mu_3-k}.
\end{align*}

We may assume 
\begin{align*}
0<\mu\equiv\mu_0-1=2\mu_1-1=\mu_2=\mu_3-1<1. 
\end{align*}

Under these settings
the Schr\"odinger operator
$H=-\triangle_g+V$ is essentially self-adjoint, and we denote 
the unique self-adjoint extension by the same symbol $H$. 
Then the equation  (\ref{7.8.15.10.0}) has the unique solution
$u(t,\cdot)=e^{-itH}u_0$ for any initial data $u_0\in {\mathcal H}$.

We now fix  $(x_0,\xi^0)\in T^*M$ and $t_0>0$.
Our main result concerns a necessary and sufficient condition
for $(x_0,\xi^0)\in \wf(e^{-it_0H}u_0)$ in terms of 
the scattering data for the classical trajectory with the initial condition $(x_0,\xi^0)$. 
We denote by $(x(t;x_0,\xi^0),\xi(t;x_0,\xi^0))
=\exp tH_p(x_0,\xi^0)$
the solution to the free Hamilton equation
on $(M,g)$:
\begin{align*}
\dot{x}=\frac{\partial p}{\partial \xi}(x,\xi),
\quad \dot{\xi}=-\frac{\partial p}{\partial x}(x,\xi),\quad
p(x,\xi)=\sum_{j,k=1}^n g^{jk}(x)\xi_j\xi_k
\end{align*}
with initial data $(x(0),\xi(0))=(x_0,\xi^0)$.
We say $(x_0,\xi_0)$ is \textit{backward nontrapping}
if $x(t;x_0,\xi_0)$ escapes from any compact set in $M$
as $t\to -\infty$,
and denote the set of all backward trapping points
by ${\mathcal T}_-$.
In the next section we show that 
if $(x_0,\xi^0)\in T^*M\setminus{\mathcal T}_-$,
then
the limits
\begin{align*}
r_-(x_0,\xi^0)&{}=\lim_{t\to-\infty}(r(t;x_0,\xi^0)-2t\rho(t;x_0,\xi^0)),&
\theta_-(x_0,\xi^0)&{}=\lim_{t\to-\infty}\theta(t;x_0,\xi^0),\\
\rho^-(x_0,\xi^0)&{}=\lim_{t\to-\infty}\rho(t;x_0,\xi^0),&
\omega^-(x_0,\xi^0)&{}=\lim_{t\to-\infty}\omega(t;x_0,\xi^0)
\end{align*}
exist, where $(r,\theta,\rho, \omega)\in T^*M_\infty$
are local coordinates
associated with the boundary decomposition
$(r,\theta)\in M_\infty$,

The limits $(r_-,\theta_-,\rho^-,\omega^-)$
are the classical scattering data at time $-\infty$.
We note $\rho^-$ is the negative square root of the total energy; 
$\theta_-$ is the asymptotic direction; 
$\omega^-$ is the impact parameter (up to a scaling factor);
and $r_-$ is the time-shift of the trajectory (times velocity).

We then set up a reference free system in the quantum mechanics.
We set
\begin{align*}
M_{\free}&{}=\R\times \partial M,\\
{\mathcal H}_{\free}&{}=(L^2(M_{\free},\sqrt{h}drd\theta),
(\cdot,\cdot)_{\free}),\\
(u,v)_{\free}
&{}=\int_{M_{\free}} u(r,\theta)\overline{v(r,\theta)}\sqrt{h(\theta)}drd\theta,\quad
h=\det (h_{jk}),
\end{align*}
where $(h_{jk})$ is the Riemannian metric on $\partial M$,
appearing in Definition~1.1. 
We now set our free Hamiltonian as 
\begin{align*}
H_0=-\frac{\partial^2}{\partial r^2},\quad
D(H_0)=\{u\in {\mathcal H}_{\free};\ H_0u\in {\mathcal H}_{\free}\}.
\end{align*}
Then $H_0$ is self-adjoint on ${\mathcal H}_{\free}$.
We choose $j\in C^\infty((1,\infty))$ such that
\begin{align*}
j(r)=\left\{
\begin{array}{ll}
1,& \mbox{if }r>2,\\
0,& \mbox{if }r<3/2,
\end{array}\right.
\end{align*}
and define $J\colon {\mathcal H}_{\free}\to {\mathcal H}$
by
\begin{align}
(Ju)(x)=\left\{
\begin{array}{ll}
j(r(x))g(x)^{-\frac{1}{4}}
h(\theta(x))^\frac{1}{4}u(r(x),\theta(x)),& \mbox{if }x\in M_\infty,\\
0,& \mbox{if }x\notin M_\infty.
\end{array}\right. 
\label{7.8.15.2.0}
\end{align}
Using the natural identification: $M_\infty\subset M_{\free}$, we may write 
\begin{align*}
(r_-,\theta_-,\rho^-,\omega^-)
=\lim_{t\to-\infty}\exp (-tH_{\rho^2})\circ \exp tH_p(x_0,\xi^0)
\in T^*M_{\free}.
\end{align*}

Now we state the main result of this paper. 

\begin{theorem}\label{7.8.20.1.45}
Assume $g$ is a scattering metric of short-range type, and let 
$H=-\triangle_g +V$ where $V$ is a smooth short-range type potential. 
Let $u_0\in {\mathcal H}$ and $t_0>0$, and
suppose $(x_0,\xi^0)\in T^*M\setminus {\mathcal T}_-$.
Then
\begin{align*}
&{}(x_0,\xi^0)\in \wf(e^{-it_0H}u_0)\\
&{}\iff (r_-(x_0,\xi^0),\theta_-(x_0,\xi^0),\rho^-(x_0,\xi^0),\omega^-(x_0,\xi^0))
\in \wf(e^{-it_0H_0}J^*u_0).
\end{align*}
\end{theorem}

\noindent
{\bf Remarks 1.}\ 
The theorem describe the wave front set of $e^{-it_0H} u_0$ in terms of the 
classical (inverse) wave operator and $e^{-it_0H_0}J^* u_0$. Since $H_0$ is 
essentially one-dimensional free Schr\"odinger operator, we can consider the time 
evolution as a linear transform: $(x,\xi)\mapsto (x+2t\xi, \xi)$ 
in the phase space, and the right hand side is easy to characterize. 
(We can also easily write the integral kernel of $e^{-itH_0}$ explicitly.)

\noindent
{\bf 2.}\ 
The theorem suggests that $e^{it H_0}J^* e^{-itH}$ is a Fourier integral operator 
corresponding to the canonical map: $(x_0,\xi_0)\mapsto (r_-,\theta_-,\rho^-,\omega^-)$. 
This is analogous to what Hassell and Wunsch proved in a different setting%
~\cite{Hassell-Wunsch-2}. We believe that this statement can be proved using our method, and 
it will be discussed in forthcoming papers. 

\noindent
{\bf 3.}\ 
As a simplest example, we consider the Euclidean space $(M,g)=(\R^n,dx^2)$,
and $H=-\triangle$ on $L^2(\R^n)$.
$(\R^n,dx^2)$ is a scattering manifold of short-range type,
identified with the interior of the half sphere
$\bar{M}=\S^{n}_+=\{w\in \R^{n+1};\ |w|=1,\ w^{n+1}\ge 0\}$ through the stereographic projection:
\begin{align*}
\sp \colon\R^n\to {\mathbb{S}}^n_+,\ 
x\mapsto
\sp(x)=
(\frac{x}{\langle x\rangle},
\frac{1}{\langle x\rangle}),\quad
\langle x\rangle =(1+|x|^2)^{1/2}.
\end{align*}
Our comparison system is 
\begin{align*}
H_0=-\frac{\partial^2}{\partial r^2}\quad \mbox{ on }L^2(\R\times \S^{n-1},
dr \cdot v_{\S^{n-1}}),
\end{align*}
where $v_{\S^{n-1}}$ is the standard
density on $\S^{n-1}$ and 
Theorem \ref{7.8.20.1.45} gives a characterization
of $\wf(e^{it\triangle}u_0)$
in terms of $\wf(e^{it\partial_r^2}r^{\frac{n-1}{4}}j(r)u_0)$.

\noindent
{\bf 4.}\ 
In  \cite{Nakamura-1} and \cite{Nakamura-3}, one of the authors studied 
the same problem for the asymptotically Euclidean case, and our result may be 
considered as a generalization of them. However, in \cite{Nakamura-1},
 $H_0=-\triangle$ is used as the free system, and hence the result in \cite{Nakamura-1} is 
not a special case of Theorem~1.2. The recent paper \cite{Nakamura-3} concerns 
the long-range case, and the comparison evolution is given by a modified free motion.

The analysis of microlocal singularities of solutions to Schr\"odinger equations with 
variable coefficients was introduced by a seminal paper by Craig, Kappeler and Strauss 
\cite{Craig-Kappeler-Strauss} in 1995. They considered Schr\"odinger  equation with 
asymptotically flat metric of short-range type on the Euclidean space, and 
showed so-called {\em microlocal smoothing property}, that is, the microlocal 
smoothness follows from the rapid decay of the initial state in a conic neighborhood 
of the asymptotic velocity (as $t\to-\infty$). This result was generalized to 
Schr\"odinger equations on manifolds with scattering metric by Wunsch~\cite{Wunsch}, 
and to long-range type perturbation (on the Euclidean space) by Nakamura~\cite{Nakamura-2}. 
We remark that both results are more precise than \cite{Craig-Kappeler-Strauss} 
in the sense that the conditions on the initial states are time-dependent. Wunsch 
used {\em quadratic scattering wave front set}\/ and Nakamura used {\em homogeneous 
wave front set} to formulate the assumption. The definitions of these notions differ 
considerably, but it was shown by Ito \cite{Ito} that they are essentially equivalent 
(up to a linear transform). These results are generalized to the microlocal analytic 
singularities by Robbiano and Zuily \cite{Robbiano-Zuily-1,Robbiano-Zuily-2,Robbiano-Zuily-3}, 
and Martinez, Nakamura and Sordoni \cite{Martinez-Nakamura-Sordoni}. 

The microlocal smoothing properties give us sufficient conditions for the microlocal 
smoothness of solutions, but they do not give us necessary conditions, which 
was addressed by Hassell and Wunsch \cite{Hassell-Wunsch-1,Hassell-Wunsch-2}. 
They considered the problem on manifold with scattering metric of short-range 
type (somewhat stronger assumptions), and gave a necessary and sufficient 
condition of the microlocal regularities of solutions in the nontrapping area. 
The same problem was considered in \cite{Nakamura-1, Nakamura-3}, and 
essentially equivalent characterization of the microlocal singularities of solutions 
was proved for Schr\"odinger equations of long-range type on the Euclidean spaces. 
An analytic analogue of these results  (short-range type on the Euclidean space) was 
recently proved by Martinez, Nakamura and Sordoni~\cite{Martinez-Nakamura-Sordoni-2}. 
The method of this paper is closely related to that of \cite{Nakamura-1,Nakamura-3}, 
though the geometric structure is quite different. In particular, in \cite{Nakamura-1}, 
the standard classical scattering theory is used to construct the scattering 
correspondence, whereas we essentially construct the classical scattering in the polar 
coordinate, which makes these results being not equivalent in the case of Euclidean space. 
We also note that our method is partially inspired by works on 
the characterization of the singularities for perturbed harmonic oscillators
(\cite{Zelditch},\cite{Kapitanski-Rodnianski-Yajima},\cite{Doi}). 
We refer \cite{Craig-Kappeler-Strauss}, \cite{Hassell-Wunsch-2} and \cite{Nakamura-3}
for other references. 

Compared to the results by Hassell and Wunsch \cite{Hassell-Wunsch-2}, our result 
is formulated using more elementary terms, and the proof seems simpler. 
Hassell and Wunsch used the multiplication operator $e^{ir^2/2t}$ to compensate 
the oscillation of the solution at $\infty$, whereas we use the one-dimensional free 
Schr\"odinger evolution $e^{-itH_0}$ for that purpose. 
Thus the correspondence of the microlocal singularities is given by the classical 
scattering relation in our paper, whereas the so-called {\em sojourn relation}\/ 
is used in \cite{Hassell-Wunsch-2}. 
Also, the scattering wave front set (due to Melrose~\cite{Melrose}) is used to characterize the 
wave front set in \cite{Hassell-Wunsch-2}, whereas we use the standard wave front set. 
We note that Hassell and Wunsch constructed a parametrix of the evolution 
operator as a Legendre distribution, whereas we have not used any parametrix 
construction but an Egorov-type theorem. 

The paper is organized as follows: In Section~2 we study the classical scattering on 
manifold with scattering metric. Section~3 is devoted to the proof of the microlocal smoothing 
property of our evolution operator. The smoothing property proved in this section is more
precise than the previous results, and it is of interest in itself (Theorem~3.2), though 
we need only its very weak version in the proof of our main result. We prove the main 
theorem in Section~4. 

\medskip
\noindent
\textbf{Acknowledgements.} 
SN wishes to thank Andrew Hassell and Jared Wunsch for valuable discussions 
and comments. He is also grateful to Kenji Yajima for his encouragement. 

\section{Classical trajectories and scattering data}
In this section we consider the asymptotic behavior of 
the
classical flow generated by the kinetic energy function
$p(x,\xi)$ on $T^*M$.
Let $(M,g)$ be a scattering manifold of short-range type,
and $(x_0,\xi^0)\in T^*M\setminus {\mathcal T}_-$,
then \[\exp tH_p (x_0,\xi^0)\in T^*M_\infty\]
for large $-t>0$. Thus we can write
\begin{align*}
\exp tH_p(x_0,\xi^0)=(r(t;x_0,\xi^0),\theta(t;x_0,\xi^0),
\rho(t;x_0,\xi^0),\omega(t;x_0,\xi^0))
\in T^*M_\infty,
\end{align*}
where
$(r(t;x_0,\xi^0),\theta(t;x_0,\xi^0),
\rho(t;x_0,\xi^0),\omega(t;x_0,\xi^0))$ satisfies
\begin{align}
\dot{r}=\frac{\partial p}{\partial \rho},\qquad
\dot{\theta}=\frac{\partial p}{\partial \omega},\qquad
\dot{\rho}=-\frac{\partial p}{\partial r},\qquad
\dot{\omega}=-\frac{\partial p}{\partial \theta}
\label{11.13.0}
\end{align}
with
\begin{align*}
&{}p(r,\theta,\rho,\omega)\\
&{}=\rho^2+\frac{1}{r^2}h^{jk}(\theta)\omega_j\omega_k
+a_0(r,\theta)\rho^2
+\frac{1}{r}a_{1}^j(r,\theta)\rho \omega_j
+\frac{1}{r^2}a_2^{jk}(r,\theta)\omega_j\omega_k.
\end{align*}
By the Cramer's formula we can
compute the inverse matrix to
\begin{align*}
(g_{jk})=
\left(
\begin{array}{cc}
1+m^{0}& r{}^tm^1\\
rm^{1} &r^2h+r^2m^2
\end{array}\right),
\end{align*}
so that we obtain
\begin{align*}
|\partial_r^k\partial_\theta^\alpha a_0(r,\theta)|
&{}\le C_{k\alpha}r^{-1-\mu-k},\\
|\partial_r^k\partial_\theta^\alpha a_1(r,\theta)|
&{}\le C_{k\alpha}r^{-(1+\mu)/2-k}\le C_{k\alpha}r^{-\mu-k},\\
|\partial_r^k\partial_\theta^\alpha a_2(r,\theta)|
&{}\le C_{k\alpha}r^{-\mu-k}.
\end{align*}
Then (\ref{11.13.0}) is written explicitly:
\begin{align}
\dot{r}={}&2\rho+2a_0\rho+\frac{1}{r}a_1^j\omega_j,
\label{4.7.14.55}\\
\dot{\theta}^j={}&\frac{2}{r^2}h^{jk}\omega_k
+\frac{1}{r}a_1^j\rho+\frac{2}{r^2}a_2^{jk}\omega_k,\label{3.19.19.19}\\
\begin{split}
\dot{\rho}={}&
\frac{2}{r^3}h^{jk}\omega_j\omega_k
-\frac{\partial a_0}{\partial r}\rho^2
+\frac{1}{r^2}a_1^j\rho\omega_j\\
&{}-\frac{1}{r}\frac{\partial a_1^j}{\partial r}\rho\omega_j
+\frac{2}{r^3}a_2^{jk}\omega_j\omega_k
-\frac{1}{r^2}\frac{\partial a_2^{jk}}{\partial r}\omega_j\omega_k,
\end{split}\label{4.7.16.0}\\
\dot{\omega}_j={}&{} 
-\frac{1}{r^2}\frac{\partial h^{kl}}{\partial \theta^j}
\omega_k\omega_l
-\frac{\partial a_0}{\partial \theta^j}\rho^2
-\frac{1}{r}\frac{\partial a_1^k}{\partial \theta^j}\rho\omega_k
-\frac{1}{r^2}\frac{\partial a_2^{kl}}{\partial \theta^j}\omega_k\omega_l.
\label{4.7.15.35}
\end{align}
\begin{lemma}\label{12.5.0}
If $(x_0,\xi^0)\in T^*M\setminus{\mathcal T}_-$, then
there is large  $C>0$ such that
\begin{align}
r(t;x_0,\xi^0) > C^{-1}|t|-C\quad \mbox{ for } t< -C.
\label{9.23.5.35}
\end{align}
The estimate holds locally uniformly in
the initial data $(x_0,\xi^0)$. In particular,
$T^*M\setminus{\mathcal T}_-$ is open.
\end{lemma}
\textit{Proof.}
First note that
\[p(r(t;x_0,\xi^0),\theta(t;x_0,\xi^0),
\rho(t;x_0,\xi^0),\omega(t;x_0,\xi^0))=
p(x_0,\xi^0)\equiv p_0\]
by the conservation of the energy.
Then, since 
\begin{align*}
p_0={}&
(\frac{1}{2}+a_0)\rho^2
+\frac{1}{r^2}(h^{jk}-\frac{1}{2}a_1^ja_1^{k}+a_2^{jk})
\omega_j\omega_k+\frac{1}{2}(\rho+\frac{1}{r}a_1^j\omega_j)^2\\
\ge{}&{}(\frac{1}{2}+a_0)\rho^2
+\frac{1}{r^2}(h^{jk}-\frac{1}{2}a_1^ja_1^{k}+a_2^{jk})
\omega_j\omega_k,
\end{align*}
we obtain
\begin{align}
|\rho(t;x_0,\xi^0)|<C,\qquad |\omega(t;x_0,\xi^0)|<Cr(t;x_0,\xi^0)\label{4.7.15.50}
\end{align}
as long as $r(t;x_0,\xi^0)$ is large. 
By direct computation we have
\begin{align*}
\ddot{r}={}&{}2\dot{\rho}+2\frac{\partial a_0}{\partial r}\dot{r}\rho
+2\frac{\partial a_0}{\partial \theta^j}\dot{\theta}^j\rho
+2a_0\dot{\rho}
-\frac{\dot{r}}{r^2}a_1^j\omega_j\\
&{}+\frac{1}{r}\frac{\partial a_1^j}{\partial r}\dot{r}\omega_j
+\frac{1}{r}\frac{\partial a_0}{\partial \theta^k}\dot{\theta}^k\omega_j
+\frac{1}{r}a_1^j\dot{\omega}_j\\
={}&{}
\frac{4}{r^3}h^{jk}\omega_j\omega_k
+O(r^{-1-\mu}),
\end{align*}
so that
\begin{align}
\frac{d^2}{dt^2}(r^2)
=2(\dot{r}^2+r\ddot{r})
=8p_0+O(r^{-\mu}).\label{9.23.5.30}
\end{align}
Noting that $T_0<0$ can be chosen so that
$r(T_0;x_0,\xi^0)$ is large and $\dot{r}(T_0;x_0,\xi^0)<0$,
we learn (\ref{9.23.5.30}) implies (\ref{9.23.5.35}).
Taking a small neighborhood of $\exp T_0H_p(x_0,\xi^0)$,
the assertion on local uniformity follows.
\hfill$\square$

\smallskip
\noindent
For any small open set 
$W\subset T^*M\setminus {\mathcal T}_-$ 
and negatively large $T_0<0$, we have 
\begin{align*}
\exp tH_p(x_0,\xi^0)
=(r(t;x_0,\xi^0),\theta(t;x_0,\xi^0),
\rho(t;x_0,\xi^0),\omega(t;x_0,\xi^0))
\in T^*M_\infty
\end{align*}
for $(x_0,\xi^0)\in W$ and $t< T_0$.
Then, identifying $M_\infty\hookrightarrow M_{\free}$,
we can define for $(x_0,\xi^0)\in W$ and $t< T_0$
\begin{align*}
&{}S_t(x_0,\xi^0)\\
&{}=\exp (-tH_{\rho^2})\circ \exp tH_p(x_0,\xi^0)\\
&{}=(r(t;x_0,\xi^0)-2t\rho(t;x_0,\xi^0),\theta(t;x_0,\xi^0),
\rho(t;x_0,\xi^0),\omega(t;x_0,\xi^0))\in T^*M_{\free}.
\end{align*}
\begin{proposition}
Let $(x_0,\xi^0)\in T^*M\setminus {\mathcal T}_-$,
then the limit
\begin{align*}
S_{-\infty}(x_0,\xi^0)
=\lim_{t\to-\infty}S_{t}(x_0,\xi^0)\in T^*M_{\free}
\end{align*}
exists.
Moreover, $S_{-\infty}\colon  T^*M\setminus {\mathcal T}_-
\to T^*M_{\free}$ is a local diffeomorphism.
\end{proposition}
\textit{Proof.}
By the equations (\ref{3.19.19.19}), (\ref{4.7.15.35})
and the estimate (\ref{4.7.15.50})
it follows that
\begin{align*}
\frac{d}{dt}(h^{jk}\omega_j\omega_k)
={}&{}\dot{\theta}^l\frac{\partial h^{jk}}{\partial \theta^l}\omega_j\omega_k
+2h^{jk}\dot{\omega}_j\omega_k\\
={}&{}\frac{2}{r^2}h^{lm}\omega_m\frac{\partial h^{jk}}{\partial \theta^l}\omega_j\omega_k
+O(r^{-1-\mu}h^{jk}\omega_j\omega_k)\\
&{}-\frac{2}{r^2}h^{jk}\frac{\partial h^{lm}}{\partial\theta^j}
\omega_l\omega_m\omega_k
+O(r^{-1-\mu}(|\omega|+h^{jk}\omega_j\omega_k))\\
={}&{}O(r^{-1-\mu}(1+h^{jk}\omega_j\omega_k) ).
\end{align*}
Then by Lemma \ref{12.5.0}
we obtain 
\begin{align*}
h^{jk}(\theta(t;x_0,\xi^0))\omega_j(t;x_0,\xi^0)\omega_k(t;x_0,\xi^0)
\le C \exp (\int r(t;x_0,\xi^0)^{-1-\mu}dt)
\le C',
\end{align*}
which holds locally uniformly in the initial data.
Thus $\omega(t;x_0,\xi^0)$ is 
bounded as $t\to -\infty$,
and we have
\begin{align*}
\frac{d}{dt}(r(t;x_0,\xi^0)-2t\rho(t;x_0,\xi^0))
=2a_0\rho+\frac{1}{r}a_1^j\omega_j-2t\dot{\rho}
=O(\langle t\rangle^{-1-\mu}),
\end{align*}
so that there exists
$r_-(x_0,\xi^0)$ such that
\begin{align*}
|r_-(x_0,\xi^0)-(r(t;x_0,\xi^0)-2t\rho(t;x_0,\xi^0))|
\le C\langle t\rangle^{-\mu}.
\end{align*}
Similarly for other variables
there exist
$\theta_-(x_0,\xi^0)$, $\rho^-(x_0,\xi^0)$ and $\omega^-(x_0,\xi^0)$
such that
\begin{align*}
|\theta_-(x_0,\xi^0)-\theta(t;x_0,\xi^0)|&{}\le C\langle t\rangle^{-\mu},\\
|\rho^-(x_0,\xi^0)-\rho(t;x_0,\xi^0)|&{}\le C\langle t\rangle^{-1-\mu},\\
|\omega^-(x_0,\xi^0)-\omega(t;x_0,\xi^0)|&{}\le C\langle t\rangle^{-\mu}.
\end{align*}
All these estimates hold locally uniformly in $(x_0,\xi^0)$,
and thus
\begin{align*}
S_{-\infty}=\lim_{t\to-\infty}S_{t}
\colon T^*M\setminus {\mathcal T}_-
\to T^*M_{\free}
\end{align*}
is locally uniform convergence and $S_{-\infty}$
is continuous. 

Next we prove the smoothness of
$S_{-\infty}$.
Since 
\begin{align}
S_{-\infty}=\exp (-tH_{\rho^2}) \circ S_{-\infty}\circ \exp tH_p
\quad \mbox{ for any }t\in \R
\label{9.23.13.0}
\end{align}
and $\exp tH_{\rho^2}$ and $\exp tH_p$
are diffeomorphisms on $T^*M_{\free}$ and $T^*M$,
respectively,
we have only to work on $T^*M_\infty\subset T^*M_{\free}$,
that is,
it suffices to show that
all the derivatives of $S_{t}$
converges as $t\to -\infty$ uniformly 
on any small open set $W\subset T^*M_\infty$ with
\begin{align*}
\exp tH_p (W)\subset  T^*M_\infty\quad \mbox{ for } t\le 0.
\end{align*}
If we put 
\begin{align*}
(R(t),\theta(t),\rho(t),\omega(t))
=\exp(-tH_{\rho^2})\circ\exp tH_p (x_0,\xi^0), \quad(x_0,\xi^0)\in W,
\end{align*}
then it satisfies
\begin{align}
\dot{R}=\frac{\partial K}{\partial \rho},\qquad
\dot{\theta}=\frac{\partial K}{\partial \omega},\qquad
\dot{\rho}=-\frac{\partial K}{\partial R},\qquad
\dot{\omega}=-\frac{\partial K}{\partial \theta},
\label{9.23.7.15}
\end{align}
where
\begin{align*}
K(t,R,\theta,\rho,\omega)={}&\rho^2-p(R+2t\rho,\theta,\rho,\theta)\\
={}&-\frac{1}{(R+2t\rho)^2}h^{jk}(\theta)\omega_j\omega_k
-a_0(R+2t\rho,\theta)\rho^2\\
&{}-\frac{a_{1}^j(R+2t\rho,\theta)}{R+2t\rho}\rho \omega_j
-\frac{a_2^{jk}(R+2t\rho,\theta)}{(R+2t\rho)^2}\omega_j\omega_k.
\end{align*}
Let us first show that all the derivatives of 
$(R(t),\theta(t),\rho(t),\omega(t))$
are bounded for $t\le 0$.
We have already shown the boundedness for
$(R(t),\theta(t),\rho(t),\omega(t))$.
Now  assume that
$(\partial^{\beta} R(t),
\partial^{\beta}\theta(t),
\partial^{\beta}\rho(t),
\partial^{\beta}\omega(t))$, $|\beta|\le l-1$
is bounded, where $\partial^{\beta}
=\partial_{(x_0,\xi^0)}^{\beta}$ for simplicity.
Then, by differentiating (\ref{9.23.7.15}), we learn
for $\alpha\in \Z^n_+$ with $|\alpha|=l$
\begin{align}
\frac{d}{dt}\partial^\alpha R(t)
=O(\langle t\rangle^{-1-\mu}
(|(\partial^\alpha R(t),
\partial^\alpha \theta(t),
\partial^\alpha \rho(t),
\partial^\alpha \omega(t))|+1)),
\label{9.22.3.30}
\end{align}
and similarly for
$\partial^\alpha \theta(t)$,
$\partial^\alpha \rho(t)$,
$\partial^\alpha \omega(t)$.
Thus, if we put
\begin{align*}
\psi(t)=|(\partial^\alpha R(t),
\partial^\alpha \theta(t),
\partial^\alpha \rho(t),
\partial^\alpha \omega(t))|,
\end{align*}
then
\begin{align*}
\psi(t)\le C+C\int_t^0\langle s\rangle^{-1-\mu}\psi(s)\, ds.
\end{align*}
The Gronwall's inequality applies and we get
\begin{align*}
\psi(t)\le C.
\end{align*}
By induction all the derivatives of 
$(R(t),\theta(t),\rho(t),\omega(t))$ are bounded,
and thus, in turn, substituting this estimate to
(\ref{9.22.3.30}) and the counterparts for
$\partial^\alpha \theta(t)$,
$\partial^\alpha \rho(t)$,
$\partial^\alpha \omega(t)$, 
we can show that
all the derivatives of $S_{t}$
converges locally uniformly as $t\to -\infty$.

Finally, in order to prove that $S_{-\infty}$
is a local diffeomorphism,
we need to estimate the first derivatives 
of $S_t$ more precisely.
Noting the boundedness of $S_{t}(W)$, $t\le 0$ 
and the relation (\ref{9.23.13.0}), we observe 
it suffices to show that
$S_{-\infty}$ is diffeomorphic on a small open set
$\tilde{W}=\exp T_0H_p(W)\subset T^*M\setminus {\mathcal T}_-$. 
For arbitrary large $R_1>0$, by taking
$T_0<0$ negatively large enough accordingly,
we may assume
\begin{align*}
S_t(\tilde{W})\subset T^*((R_1,\infty)\times \partial M)
\setminus {\mathcal T}_-.
\end{align*}
Thus similarly to (\ref{9.22.3.30})
we obtain for $(R_0,\theta_0,\rho^0,\omega^0)\in \tilde{W}$
\begin{align*}
&{}\left\lVert \frac{d}{dt}\left(\frac{
\partial (R(t),\theta(t),\rho(t),\omega(t))}{
\partial (R_0,\theta_0,\rho^0,\omega^0)}- 1\right)\right\rVert\\
&{}\le O(\langle R_1+2t\rho \rangle^{-1-\mu})
\left[
\left\lVert \frac{
\partial (R(t),\theta(t),\rho(t),\omega(t))}{
\partial (R_0,\theta_0,\rho^0,\omega^0)}- 1\right\rVert
+1\right],
\end{align*}
and 
\begin{align*}
\frac{
\partial (R(0),\theta(0),\rho(0),\omega(0))}{
\partial (R_0,\theta_0,\rho^0,\omega^0)}= 1,
\end{align*}
so that, if $R_1$ is large,
\begin{align*}
\left\lVert \frac{
\partial (R(t),\theta(t),\rho(t),\omega(t))}{
\partial (R_0,\theta_0,\rho^0,\omega^0)}- 1\right\rVert
<\frac{1}{2}\quad \mbox{ for }t\le 0.
\end{align*}
It follows that the Jacobian of
$S_{-\infty}$ is invertible on $\tilde{W}$, and thus
by (\ref{9.23.13.0})
$S_{-\infty}$ is a local diffeomorphism.
\hfill
$\square$

\section{Microlocal smoothing estimates}
Note that, if a boundary decomposition $M_\infty\cong(1,\infty)
\times \partial M$ is fixed, 
the scaling 
\[T^*M_\infty\to T^*M_\infty,\ (r,\theta,\rho,\omega)\mapsto
(h^{-1}r,\theta,h^{-1}\rho,h^{-1}\omega), \quad h\in (0,1]\]
is well-defined.
Thus, if
$a\in C_0^\infty(T^*M_\infty)$,
then $a(hr,\theta,h\rho,h\omega)$,
$h\in (0,1]$ is well-defined.
Choose $\psi\in C^\infty(M)$
so that it does not depend on $r$ near $\partial M$,
equals $1$ on
\begin{align*}
\{(h^{-1}r,\theta)\in M_\infty;\ 
h\in (0,1],\ (r,\theta)\in \pi (\supp a)\},
\end{align*}
and has support sufficiently near the same set.
Here $\pi \colon T^*M\to M$ is the projection to
the base space.
If $\supp a$ is sufficiently small, 
we can choose $\psi$ and a coordinate neighborhood
$U\subset \partial M$
such that $\supp \psi\subset (1,\infty)\times U$.
Using the chart diffeomorphism
$\kappa \colon (1,\infty)\times U\stackrel{\cong}{\to}
(1,\infty)\times \tilde{U}\subset \R^n$, 
we have for $v\in {\mathcal H}$
\begin{align*}
\kappa^*[
(\kappa_*g)^{-\frac{1}{4}}
(\kappa_*\psi)
(\kappa_*a)^{\w}(hr,\theta,hD_r,hD_\theta)
(\kappa_*\psi)
(\kappa_*g)^{\frac{1}{4}}
(\kappa_*v)]\in {\mathcal H},
\end{align*}
where 
$\kappa^*$ and $\kappa_*=(\kappa^*)^{-1}$
are
the associated pull-back and the push-forward, respectively.
In the following argument there always appears 
only one coordinate chart.
Hence $(1,\infty)\times U$ may be identified with 
$(1,\infty)\times \tilde{U}\subset\R^n$, and
$\kappa^*$ and $\kappa_*$ will not be written explicitly,
if there is no confusion.
\begin{definition}
Let $(r_0,\theta_0,\rho^0,\omega^0)\in T^*M_\infty$
and $v\in {\mathcal S}'(M)$.
We write 
\begin{align*}
(r_0,\theta_0,\rho^0,\omega^0)\notin \wf^{\rh}(v)
\end{align*}
if there are $a\in C_0^\infty(T^*M_\infty)$
with $a(r_0,\theta_0,\rho^0,\omega^0)\neq 0$
and $\psi \in C^\infty(M)$ as above
such that 
\begin{gather*}
\lVert g^{-\frac{1}{4}}\psi
a^{\w}(hr,\theta,hD_r,hD_\theta)
\psi g^{\frac{1}{4}}v\rVert_{{\mathcal H}}=O(h^\infty)
\quad \mbox{ as }h\downarrow 0,
\end{gather*}
(where we have omitted $\kappa$).
We call $\wf^{\rh}(v)$ the \textit{radially homogeneous
wave front set} of $v$.
\end{definition}

\noindent
{\bf Remark}\ 
Recall the semiclassical characterization of the 
(usual) wave front set, i.e., 
$(x_0,\xi^0)\notin \wf(v)$ if and only
if there are coordinate neighborhood $V$ around
$x_0$, $a\in C^\infty_0(T^*V)$
with $a(x_0,\xi^0)\neq 0$
and $\psi\in C^\infty_0(V)$ with $\psi=1$ on $\pi(\supp a)$
such that
\begin{align*}
\lVert g^{-\frac{1}{4}}\psi a^{\w}(x,hD_x)
\psi g^{\frac{1}{4}}v\rVert_{{\mathcal H}}=O(h^\infty)
\quad \mbox{ as }h\downarrow 0.
\end{align*}
(See also the \textit{frequency set}, Definition \ref{9.24.13.40}.)
The radially homogeneous wave front set differs
from the wave front set only in the parameter
in front of $r$ in the polar coordinates.

The \lq\lq homogeneous'' wave front set 
corresponds to a somewhat different set corresponding
to the operators of the form
\begin{align*}
a^{\w}(hr,\theta,hD_r,h^2D_\theta).
\end{align*}
The homogeneous wave front set is used
for the long-range case in \cite{Nakamura-2} and \cite{Ito}.

\begin{theorem}\label{13.5.20}
Let $u_0\in {\mathcal H}$, and $t_0>0$.
Suppose $(x_0,\xi^0)\in T^*M\setminus {\mathcal T}_-$, and
denote its scattering data at $t=-\infty$
by $(r_-,\theta_-,\rho^-,\omega^-)$.
If 
\begin{align*}
(-2t_0\rho^-,\theta_-,\rho^-,\omega^-)\notin \wf^{\rh}(u_0),
\end{align*}
then we have 
\begin{align*}
(x_0,\xi^0)\notin \wf(e^{-it_0H}u_0).
\end{align*}
\end{theorem}
Theorem \ref{13.5.20}
is a refinement of the results of \cite{Nakamura-2}
and \cite{Ito} in the short-range case,
and also they are refinements of the 
microlocal smoothing property
of Craig-Kappeler-Strauss \cite{Craig-Kappeler-Strauss} and Wunsch~\cite{Wunsch}.
While \cite{Craig-Kappeler-Strauss} requires
the initial data to decay 
in some cone in the configuration space
in its assumption,
the homogeneous or radially homogeneous wave front set
capture the decaying property in cones in the phase space,
which microlocally weaken the assumption of \cite{Craig-Kappeler-Strauss}.
\begin{corollary}
Let $J$ be as in (\ref{7.8.15.2.0}).
Then
\begin{align*}
\wf(e^{-itH}u_0)\setminus{\mathcal T}_-
=\wf(e^{-itH}JJ^*u_0)\setminus{\mathcal T}_-\quad 
t>0.
\end{align*}
\end{corollary}
\textit{Proof.}
Note that the operator $JJ^*$ is the
multiplication operator by $j(x)^2$.
Then we decompose
\begin{align*}
e^{-itH}u_0=e^{-itH}JJ^*u_0+e^{-itH}(1-JJ^*)u_0.
\end{align*}
It is easy to see $\wf^{\rh}((1-JJ^*)u_0)=\varnothing$,
so that
\[\wf(e^{-itH}(1-JJ^*)u_0)\setminus {\mathcal T}_-=\varnothing
\quad t>0.\]
Then the assertion follows.
\hfill$\square$

\smallskip
\noindent
The proof 
of Theorem \ref{13.5.20}
is analogous to \cite{Nakamura-2},
and rather long.
We first prove some lemmas needed later.
\begin{definition}\label{9.24.13.40}
Let $v_h\in {\mathcal H}$ be a vector dependent on $h\in (0,1]$ with 
$\lVert v_h\rVert \le 1$, and $(x_0,\xi^0)\in T^*M$.
We say that $v_h$ is \textit{microlocally infinitely small} near
$(x_0,\xi^0)$, if there exist a chart
$\kappa\colon V\to \tilde{V}\subset \R^n$ around $x_0$,
$a\in C^\infty_0(T^*V)$ with $a(x_0,\xi^0)\neq 0$,
and a cutoff function $\psi\in C_0^\infty(V)$
with $\psi=1$ on $\supp a$ such that
\begin{align*}
\lVert g^{-\frac{1}{4}}\psi a^{\w}(x,hD_x)
\psi g^{\frac{1}{4}}
v_h\rVert_{{\mathcal H}}=O(h^\infty).
\end{align*}
The \textit{frequency set} $\fs (v_h)$ is the compliment 
in $T^*M$ of
such $(x_0,\xi^0)$'s.
\end{definition}
\begin{lemma}\label{7.8.20.9.15}
Let $(x_0,\xi^0)\in T^*M$ and $v\in {\mathcal H}$.
Then for any $T_0\in \R$
\begin{align}
(x_0,\xi^0)\in \wf(v)\iff
\exp T_0H_p(x_0,\xi^0)\in \fs(e^{-ihT_0H}v).
\label{9.24.16.30}
\end{align}
\end{lemma}
\textit{Proof.} 
Let $\kappa\colon V\to \tilde{V}\subset\R^n$ be any chart
around $(x_0,\xi^0)$. We assume $T_0<0$ and
$\exp tH_p(x_0,\xi^0)\in V$ for $t\in [T_0,0]$.
Suppose $a\in C^\infty_0(T^*V)$
has sufficiently small support near $(x_0,\xi^0)$,
and put
\begin{align*}
b_0(x,\xi;t)=a\circ \exp (-tH_p)(x,\xi)\quad
(x,\xi,t)\in T^*V\times [T_0,0].
\end{align*}
Note that $b_0$ is the solution to
the first transport equation:
\begin{align*}
\frac{\partial }{\partial t}b_0+\{p,b_0\}=0,\quad
b_0(\cdot,\cdot;0)=a,
\end{align*}
where $\{\cdot,\cdot\}$ is the Poisson bracket.
Since $\supp a$ is small, we may assume
\begin{align*}
\bigcup_{t\in [T_0,0]}
\supp b_0(\cdot,\cdot; t)\subset T^*V.
\end{align*}
Take $\psi_0 \in C_0^\infty(V)$ with $\psi_0=1$ 
sufficiently near
$\cup_{t\in [T_0,0]}\pi(\supp b(\cdot,\cdot; t))$,
and define
\begin{align*}
B_0(t,h)=[g(x)^{-\frac{1}{4}}\psi_0(x) b_0^{\w}(x,hD_x;t)
\psi_0(x) g(x)^{\frac{1}{4}}]^2,
\end{align*}
then, by the standard argument in the
pseudodifferential calculus (see, e.g. \cite{Martinez}),
\begin{align*}
\frac{d}{dt}B_0(t,h)+ih[H,B_0(t,h)]
=g(x)^{-\frac{1}{4}}\psi_1(x) r_0^{\w}(x,hD_x;t,h)
\psi_1(x) g(x)^{\frac{1}{4}},
\end{align*}
where $r_0\in S(h)$ is supported in 
$\supp b_0$ modulo $S(h^\infty)$,
and $\psi_1\in C_0^\infty(V)$
is any function 
with $\psi_1=1$ on $\supp \psi_0$.
Here we have put $S(h)=S(h,dx^2+d\xi^2)$
in H\"ormander's notation in the chart we are concerned,
and the dependence on the parameters $t,h$ are
supposed to be uniform.
(This symbol class
is not invariant under the coordinates change,
but we do not have to take care of this fact,
since we have fixed one chart.)
Next we decompose
$r_0=r_0'+r_0''$ so that
\begin{align*}
\supp r_0'\subset\supp b_0,\quad r_0''\in S(h^\infty),
\end{align*}
and consider the second transport equation:
\begin{align*}
\frac{\partial }{\partial t}b_1+\{p,b_1\}=r'_0,\quad
b_1(\cdot,\cdot;0,h)=0.
\end{align*}
The solution is written as
\begin{align*}
b_1(x,\xi;t,h)
=\int^{t}_0 r_0'(
\exp (s-t)H_p(x,\xi); s,h)\,ds
\in S(h).
\end{align*}
Let $\supp \psi_1$ be sufficiently near $\supp \psi_0$
and put 
\begin{align*}
B_1(t,h)=g(x)^{-\frac{1}{4}}\psi_1(x)
b_1^{\w}(x,hD_x;t)
\psi_1(x) g(x)^{\frac{1}{4}},
\end{align*}
then we obtain
\begin{align*}
&{}\frac{d}{dt}(B_0(t,h)+B_1(t,h))+ih[H,B_0(t,h)+B_1(t,h)]\\
&{}=g^{-\frac{1}{4}}(x)\psi_2(x) r_1^{\w}(x,hD_x;t,h)
\psi_2(x) g(x)^{\frac{1}{4}}.
\end{align*}
Here, similarly to the above,
$r_1\in S(h^2)$ is supported in 
$\supp b_0$ modulo $S(h^\infty)$,
and $\psi_2\in C_0^\infty(V)$
is any function 
with $\psi_2=1$ on $\supp \psi_1$.
Iterating this procedure and putting
\begin{align*}
b(x,\xi;t,h)\sim \sum_{j=0}^\infty b_j(x,\xi;t,h),&&
B(t,h)=g(x)^{-\frac{1}{4}}\psi(x)
b_1^{\w}(x,hD_x;t)
\psi(x) g(x)^{\frac{1}{4}},
\end{align*}
where $\psi\in C_0^\infty(V)$
is any function 
with $\psi=1$ on $\cup_{j=0}^\infty\supp \psi_j$,
we obtain
\begin{gather*}
\frac{d}{dt}B(t,h)+ih[H,B(t,h)]
=O(h^\infty)\quad \mbox{ in }{\mathcal L}({\mathcal H}),\\
B(0,h)=[g(x)^{-\frac{1}{4}}\psi_0(x)a^{\w}(x,hD_x)
\psi_0(x) g(x)^{\frac{1}{4}}]^2.
\end{gather*}
Also note that $B(T_0,h)$ is an $h$-pseudodifferential
operator
with principal symbol $a\circ \exp (-T_0H_p)(x,\xi)^2$.
Then 
\begin{align*}
&{}(e^{-ihT_0H}v,B(T_0,h)e^{-ihT_0H}v)_{{\mathcal H}}
-(v,B(0,h)v)_{{\mathcal H}}\\
&{}=\int_0^{T_0}\frac{d}{dt}
(e^{-ihtH}v,B(t,h)e^{-ihtH}v)_{{\mathcal H}}\,dt
=O(h^\infty),
\end{align*}
and thus (\ref{9.24.16.30}) holds.

For general $T_0\in\R$, divide the geodesic
into a finite number of small segments each of which
is contained in some coordinate neighborhood,
and apply the above argument on each chart,
then the lemma follows.
\hfill$\square$

\smallskip
\noindent
Let $u_0\in {\mathcal H}$,
$(x_0,\xi^0)\in T^*M\setminus {\mathcal T}_-$ and 
$t_0>0$  be as in the assumption of Theorem \ref{13.5.20}.
By Lemma \ref{7.8.20.9.15}
it suffices to show
\begin{align*}
\exp T_0H_p (x_0,\xi^0)\notin
\fs(e^{-i(t_0+hT_0)H}u_0)
\end{align*}
for some $T_0\in \R$.
Let us choose sufficiently large negative $T_0\ll 0$
so that
we can work in a fixed coordinate neighborhood
$(1,\infty)\times U\subset M_\infty$, $U\subset \partial M$
near the trajectory:
\[\exp tH_p (x_0,\xi^0)
\in (1,\infty)\times U \quad\mbox{for } t\le T_0.\]
We put 
\begin{align*}
(r(t),\theta(t),\rho(t),\omega(t))=
\exp (t+T_0)H_p(x_0,\xi^0).
\end{align*}
Fix small $\delta'>0$ and large $C>0$,
specified later.
Take any $\delta_0\in (\delta'/2,\delta')$
and define 
\[\varphi\colon (-\infty,0]\times T^*M\rightarrow \R\]
by
\begin{align*}
\begin{split}
\varphi(t,r,\theta,\rho,\omega)
&{}= \chi\left(\frac{|r-r(t)|}{5\delta_0| t+T_0|}\right)
\chi\left(\frac{|\theta-\theta(t)|}{\delta_0-C|t+T_0|^{-\mu}}\right)\\
&{}\phantom{{}={}}\chi\left(\frac{|\rho-\rho(t)|}{\delta_0-C|t+T_0|^{-\mu-1}}\right)
\chi\left(\frac{|\omega-\omega(t)|}{\delta_0-C |t+T_0|^{-\mu}}\right)\\
&=: \chi_1\chi_2\chi_3\chi_4,
\end{split}
\end{align*}
where
\begin{align*}
\chi\in C^\infty([0,+\infty)),&&
\chi(\tau)=\left\{
\begin{array}{ll}
1,& \mbox{ if } \tau \le 1,\\
0,& \mbox{ if } \tau \ge 2,
\end{array}\right. && 
\frac{d\chi}{d\tau}\le 0.
\end{align*}
\begin{lemma}\label{13.5.5}
Let small $\delta'$ and large $C$ be given,
and take $T_0$ large enough accordingly.
Then for any $\delta_0\in (\delta'/2,\delta')$
$\varphi$ is well-defined and satisfies the following:
\begin{enumerate}
\item 
\[\varphi(t,r,\theta,\rho,\omega)\ge 0 \quad \mbox{on}\quad
(-\infty,0]\times T^*M,\]
and
\[\varphi(t,r(t),\theta(t),\rho(t),\omega(t))=1
\quad \mbox{for all}\quad t\le 0.\]
\item \label{9.10.9.50} The inequality
\[\frac{D}{Dt}\varphi(t,r,\theta,\rho,\omega)\le 0
\quad\mbox{on}\quad (-\infty,0]\times T^*M\]
holds, where
$\frac{D}{Dt}$ is the \textit{Lagrange derivative} defined by
\[\frac{D}{Dt}=\frac{\partial}{\partial t}+\{p,\cdot\}
=\frac{\partial}{\partial t}
+\frac{\partial p}{\partial \rho}\frac{\partial}{\partial r}
+\frac{\partial p}{\partial \omega}\frac{\partial}{\partial \theta}
-\frac{\partial p}{\partial r}\frac{\partial}{\partial \rho}
-\frac{\partial p}{\partial \theta}\frac{\partial}{\partial \omega}.\]
\item The inequality
\begin{align*}
|\partial_r^k\partial_\theta^\alpha
\partial_\rho^l\partial_\omega^\beta
\partial_t^n\varphi|
\le C_{k\alpha l\beta n}
\langle t\rangle^{-n-k}
\end{align*}
holds,
that is,
\[\partial_t^n
\varphi\in S(\langle t\rangle^{-n},\langle t\rangle^{-2}dr^2+d\theta^2
+d\rho^2+d\omega^2)\quad \mbox{
uniformly in }t\le 0.\]
\end{enumerate}
\end{lemma}
\textit{Proof.}
\textit{1.} is obvious. 

\smallskip
\noindent
\textit{2.}
Note that
\begin{align*}
\frac{D}{Dt}\varphi
= \frac{D\chi_1}{Dt}\chi_2\chi_3\chi_4
+\chi_1\frac{D\chi_2}{Dt}\chi_3\chi_4
+\chi_1\chi_2\frac{D\chi_3}{Dt}\chi_4
+\chi_1\chi_2\chi_3\frac{D\chi_4}{Dt}.
\end{align*}
Now let us compute the differentiations.
We first get
\begin{align*}
\frac{D\chi_1}{Dt}
=\frac{1}{5\delta_0|t+T_0|}\left[
\frac{|r-r(t)|}{|t+T_0|}
+\frac{r-r(t)}{|r-r(t)|}
\left(\frac{\partial p}{\partial \rho}
-\frac{\partial p}{\partial \rho}(t)
\right)\right]
\chi'\left(\frac{|r-r(t)|}{5\delta_0 |t+T_0|}\right).
\end{align*}
On $\supp (\frac{D\chi_1}{Dt}\chi_2\chi_3\chi_4)$
we have
\begin{align*}
1\le\frac{|r-r(t)|}{5\delta_0 |t+T_0|}\le 2
\end{align*}
and
\begin{align*}
\left|\frac{\partial p}{\partial \rho}
-\frac{\partial p}{\partial \rho}(t)\right| 
\le 4\delta_0+O(|t+T_0|^{-1-\mu}).
\end{align*}
Thus, taking $T_0$ larger if necessary,
we obtain 
\begin{align*}
\frac{D\chi_1}{Dt}\chi_2\chi_3\chi_4
\le \frac{1}{5\delta_0|t+T_0|}(\delta_0
+O(|t+T_0|^{-1-\mu}))
\chi'\chi_2\chi_3\chi_4\le 0.
\end{align*}
Similarly, by direct computations,
\begin{align*}
\frac{D\chi_2}{Dt}
={}& \frac{1}{\delta_0-C|t+T_0|^{-\mu}}\\
&{}\left[
C\mu
|t+T_0|^{-1-\mu}\frac{|\theta-\theta(t)|}{\delta_0-C|t+T_0|^{-\mu}}
+\frac{\theta-\theta(t)}{|\theta-\theta(t)|}
\left(\frac{\partial p}{\partial \omega}
-\frac{\partial p}{\partial \omega}(t)\right)
\right]
\chi',
\\
\frac{D\chi_3}{Dt}
={}&  \frac{1}{\delta_0-C|t+T_0|^{-1-\mu}}\\
&{}\left[
C(1+\mu)
|t+T_0|^{-2-\mu}\frac{|\rho-\rho(t)|}{\delta_0-C|t+T_0|^{-1-\mu}}
-\frac{\rho-\rho(t)}{|\rho-\rho(t)|}
\left(\frac{\partial p}{\partial r}
-\frac{\partial p}{\partial r}(t)\right)
\right]
\chi',
\\
\frac{D\chi_4}{Dt}
={}&  \frac{1}{\delta_0-C|t+T_0|^{-\mu}}\\
&{}\left[
C\mu
|t+T_0|^{-1-\mu}\frac{|\omega-\omega(t)|}{\delta_0-C|t+T_0|^{-\mu}}
-\frac{\omega-\omega(t)}{|\omega-\omega(t)|}
\left(\frac{\partial p}{\partial \theta}
-\frac{\partial p}{\partial \theta}(t)\right)
\right]
\chi',
\end{align*}
and, on $\supp (\chi_1\chi_2\chi_3\chi_4)$,
\begin{align*}
\left|\frac{\partial p}{\partial \omega}
-\frac{\partial p}{\partial \omega}(t)\right| 
&{}\le O(|t+T_0|^{-1-\mu}),\\
\left|\frac{\partial p}{\partial r}
-\frac{\partial p}{\partial r}(t)\right| 
&{}\le O(|t+T_0|^{-2-\mu}),\\
\left|\frac{\partial p}{\partial \theta}
-\frac{\partial p}{\partial \theta}(t)\right| 
&{}\le O(|t+T_0|^{-1-\mu}).
\end{align*}
Thus we have
\begin{align*}
\chi_1\frac{D\chi_2}{Dt}\chi_3\chi_4
&\le \frac{1}{\delta_0-C|t+T_0|^{-\mu}}\left[
C\mu|t+T_0|^{-1-\mu}+O
(|t+T_0|^{-1-\mu})\right]
\chi_1\chi'
\chi_3\chi_4,\\
\chi_1\chi_2\frac{D\chi_3}{Dt}\chi_4
&\le \frac{1}{\delta_0-C|t+T_0|^{-1-\mu}}\left[
C(1+\mu) |t+T_0|^{-2-\mu}
+O(|t+T_0|^{-2-\mu})\right]
\chi_1\chi_2
\chi'\chi_4,
\\
\chi_1\chi_2\chi_3\frac{D\chi_4}{Dt}
&\le 
\frac{1}{\delta_0-C|t+T_0|^{-\mu}}\left[
C\mu|t+T_0|^{-1-\mu}+
O(|t+T_0|^{-1-\mu})\right]
\chi_1\chi_2\chi_3\chi'.
\end{align*}
Since $C$ is large enough, we obtain the 
assertion.

\smallskip
\noindent
\textit{3.}
The estimates follow from direct computations.
\hfill$\square$

\bigskip
\noindent
We put
\[b_0(r,\theta,\rho,\omega;t,h)
=\varphi(h^{-1}t,r,\theta,\rho,\omega)
\in S(1,\langle h^{-1}t\rangle^{-2}dr^2+d\theta^2
+d\rho^2+d\omega^2),\]
and we restrict the parameter $t\in(-\infty,0]$ to the interval
$[-t_0,0]$.
Let $\psi_0\in C^\infty(M_\infty)$
be such that it does not depend on $r$ near $\partial M$,
equals $1$ on 
\begin{align*}
\bigcup
\{\pi (\supp b_0(\cdot,\cdot,\cdot,\cdot;t,h));\ 
t\in [-t_0,0], h\in (0,1]\}
\end{align*}
and has support sufficiently near the same set.
We define
\begin{align*}
B_0(t,h)=[g(r,\theta)^{-\frac{1}{4}}
\psi_0(r,\theta)
b_{0}^{\w}(r,\theta,hD_r,hD_\theta;t,h)
\psi_0(r,\theta)
g(r,\theta)^{\frac{1}{4}}]^2
\quad \mbox{ on } {\mathcal H}.
\end{align*}
\begin{lemma}\label{12.2.30}
There exists $r_0\in S(1,
\langle h^{-1}t\rangle^{-2}dr^2+d\theta^2
+d\rho^2+d\omega^2)$
such that
\begin{align*}
&{}\frac{d}{dt}B_0(t,h)+i[H,B_0(t,h)]\\
&{}\le g(r,\theta)^{-\frac{1}{4}}
\psi_1(r,\theta)
r_{0}^{\w}(r,\theta,hD_r,hD_\theta;t,h)
\psi_1(r,\theta)
g(r,\theta)^{\frac{1}{4}},
\end{align*}
and that $r_0$ has support in
$\supp b_0$ modulo $S(h^\infty)$.
Here $\psi_1\in C^\infty(M_\infty)$
is any function that
does not depend on $r$ near $\partial M$,
equals $1$ on $\supp \psi_0$, and 
has support sufficiently near the same set.
\end{lemma}
\textit{Proof.}
We compute the principal symbol of
\[\frac{d}{dt}B_0(t,h)+i[H,B_0(t,h)]
=\frac{d}{dt}B_0(t,h)
-i[\triangle_g,B_0(t,h)]
+i[V,B_0(t,h)],\]
and apply the sharp G\r{a}rding inequality.
Let $\psi_1\in C^\infty(M_\infty)$ be as in the assertion.
By the standard argument in the pseudodifferential calculus
we can write
\begin{align*}
&{}\frac{d}{dt}B_0(t,h)
-i[\triangle_g,B_0(t,h)]\\
&{}=
g^{-\frac{1}{4}}
\psi_1
[
2b_0
(\partial_t b_0
+h^{-1}\{p,b_0\})+r_{0,1}
]^{\w}(r,\theta,hD_r,hD_\theta;t,h)\psi_1
g^{\frac{1}{4}},
\end{align*} 
where the remainder term $r_{0,1}\in S(1,
\langle h^{-1}t\rangle^{-2}dr^2+d\theta^2
+d\rho^2+d\omega^2)$
has the support property
\[\supp r_{0,1}\subset \supp b_0 \mod
S(h^\infty).\]
It is also easy to see
\begin{align*}
i[V,B_0(t,h)]
&{}=
g(r,\theta)^{-\frac{1}{4}}
\psi_1(r,\theta)
r_{0,2}^{\w}(r,\theta,hD_r,hD_\theta;t,h)
\psi_1(r,\theta)
g(r,\theta)^{\frac{1}{4}},\\
r_{0,2}&{}\in S(h^{\mu},
\langle h^{-1}t\rangle^{-2}dr^2+d\theta^2
+d\rho^2+d\omega^2),\\
\supp r_{0,2}&{}\subset \supp b_0 \mod
S(h^\infty).
\end{align*}
Thus, by Lemma \ref{13.5.5},
the principal part is $2b_0
(\partial_t b_0
+h^{-1}\{p,b_0\})$, and it satisfies
\begin{align*}
2b_0
(\partial_t b_0
+h^{-1}\{p,b_0\})\in 
S(h^{-1},
\langle h^{-1}t\rangle^{-2}dr^2+d\theta^2
+d\rho^2+d\omega^2),\ {}\le 0.
\end{align*}
By the sharp G\r{a}rding inequality
$r_{0,3}\in S(1,
\langle h^{-1}t\rangle^{-2}dr^2+d\theta^2
+d\rho^2+d\omega^2)$
is found such that
\begin{gather*}
\supp r_{0,3}\subset \supp b_0 \mod
S(h^\infty)
\end{gather*}
and that
\begin{align*}
2b_0(\partial_t b_0
+h^{-1}\{p,b_0\})^{\w}
(r,\theta,hD_r,hD_\theta;t,h)
\le r_{0,3}^{\w}(r,\theta,hD_r,hD_\theta;t,h).
\end{align*}
Putting the remainder terms $r_{0,j}$, $j=1,2$
together, we obtain
\begin{align*}
&{}\frac{d}{dt}B_0(t,h)+i[H,B_0(t,h)]\\
&{}\le 
g(r,\theta)^{-\frac{1}{4}}
\psi_1(r,\theta)
r_{0}^{\w}(r,\theta,hD_r,hD_\theta;t,h)
\psi_1(r,\theta)
g(r,\theta)^{\frac{1}{4}},
\end{align*}
where $r_{0}=r_{0,1}+r_{0,2}+r_{0,3}$.
\hfill$\square$

\smallskip
Let $\delta'$, $C$ and $T_0$ be as in Lemma
\ref{13.5.5}.
Take an increasing sequence
\[\frac{\delta'}{2}<\delta_0<\delta_1<\delta_2<\cdots<\delta<\delta',\]
and set
\begin{align*}
b_j(r,\theta,\rho,\omega;t,h)
={}& \chi\left(\frac{|r-r(h^{-1}t)|}{5\delta_j| h^{-1}t+T_0|}\right)
\chi\left(\frac{|\theta-\theta(h^{-1}t)|}{\delta_j-C|h^{-1}t+T_0|^{-\mu}}\right)\\
&{}\cdot\chi\left(\frac{|\rho-\rho(h^{-1}t)|}{\delta_j-C|h^{-1}t+T_0|^{-1-\mu}}\right)
\chi\left(\frac{|\omega-\omega(h^{-1}t)|}{\delta_j-C |h^{-1}t+T_0|^{-\mu}}\right).
\end{align*}
Let $\psi_j\in C^\infty(M_\infty)$
be such that it does not depend on $r$ near $\partial M$
and that
\begin{gather*}
\psi_j=1\quad \mbox{ on }\bigcup
\{\pi (\supp b_j(\cdot,\cdot,\cdot,\cdot;t,h));\ 
t\in [-t_0,0], h\in (0,1]\},\\
\supp \psi_j\subset
\bigcup
\{\pi (\supp b_{j+1}(\cdot,\cdot,\cdot,\cdot;t,h));\ 
t\in [-t_0,0], h\in (0,1]\}.
\end{gather*}
We set for $j=1,2,\dots$
\begin{gather*}
B_j(t,h)= g(r,\theta)^{-\frac{1}{4}}
\psi_j(r,\theta)
C_jh^{\mu (j-1)}|t|b_j^{\w}(r,\theta,hD_r,hD_\theta;t,h)
\psi_j(r,\theta)
g(r,\theta)^{\frac{1}{4}}.
\end{gather*}
The constants $C_j$ will be determined inductively.
Since $b_1$ is bounded from below by a positive constant on 
$\supp b_0$,
there is large $C_1>0$ such that
\[r_0'(r,\theta,\rho,\omega;t,h)
\le C_1b_1(r,\theta,\rho,\omega;t,h),\]
where we have decomposed $r_0=r_0'+r_0''$ so that
\begin{align*}
\supp r_0'\subset\supp b_0,\quad r_0''\in S(
h^\infty).
\end{align*}
Now let us consider the operator
$\frac{d}{dt}B_1+i[H,B_1]$,
and iterate the argument in the proof of Lemma \ref{12.2.30}.
Noting that
\begin{align*}
\frac{\partial p}{\partial \theta}&{}\in
S(\langle r\rangle^{-1-\mu}\langle \rho;\omega\rangle^2,
\langle r\rangle^{-2}dr^2+d\theta^2
+\langle \rho;\omega\rangle^{-2}d\rho^2
+\langle \rho;\omega\rangle^{-2}d\omega^2),\\
\frac{\partial p}{\partial \omega}&{}\in
S(\langle r\rangle^{-1-\mu}\langle \rho;\omega\rangle,
\langle r\rangle^{-2}dr^2+d\theta^2
+\langle \rho;\omega\rangle^{-2}d\rho^2
+\langle \rho;\omega\rangle^{-2}d\omega^2)
\end{align*}
and $r=O(\langle h^{-1}t\rangle)$ on $\supp b_1$, 
we can write
as operators on $L^2(\R^n)$
through the chart diffeomorphism
\begin{align*}
&{}g(r,\theta)^{\frac{1}{4}}\Bigl(
\frac{d}{dt}B_1+i[H,B_1]\Bigr)
g(r,\theta)^{-\frac{1}{4}}\\
&{}=
\Bigl[C_1|t|\frac{D b_1}{D t}-C_1b_1
+r_{1,1}
\Bigr]^{\w}(r,\theta,hD_r,hD_\theta;t,h),
\end{align*}
where
\begin{align*}
C_1|t|\frac{D b_1}{D t}-C_1b_1&{}\in S(1,\langle h^{-1}t\rangle^{-2}dr^2+d\theta^2
+d\rho^2+d\omega^2),\\
r_{1,1}&{}\in S(h^{\mu},\langle h^{-1}t\rangle^{-2}dr^2+d\theta^2
+d\rho^2+d\omega^2).\\
\supp r_{1,1}&{}\subset \supp b_1\quad \mod S(h^\infty).
\end{align*}
Since
\begin{align*}
C_1|t|\frac{D b_1}{D t}-(
C_1b_1-r_0')
\le 0,
\end{align*}
we can find by the sharp G\r{a}rding inequality
\begin{gather*}
r_{1,2}\in S(h,\langle h^{-1}t\rangle^{-2}dr^2
+d\theta^2
+d\rho^2+d\omega^2),\\
\supp r_{1,2}\subset \supp b_1\quad \mod
S( h^\infty)
\end{gather*}
such that
\begin{align*}
\Bigl[C_1|t|\frac{D b_1}{D t}-(
C_1b_1-r_0')\Bigr]^{\w}(r,\theta,hD_r,hD_\theta;t,h)
\le r_{1,2}^{\w}(r,\theta,hD_r,hD_\theta;t,h).
\end{align*}
Also here
the inequality is in the operator sense on $L^2(\R^n)$.
Combining with the remainder terms $r_{0}''$ and $r_{1,1}$,
we obtain
\begin{align*}
&{}
\frac{d}{dt}B_1(t,h)+i[H,B_1(t,h)]
\\
&{}\le g(r,\theta)^{-\frac{1}{4}}\psi_1(r,\theta)
r_1^{\w}(r,\theta,hD_r,hD_\theta;t,h)
\psi_1(r,\theta)g(r,\theta)^{\frac{1}{4}}\\
&{}\phantom{{}={}}{}-g(r,\theta)^{-\frac{1}{4}}
\psi_0(r,\theta)
r_0^{\w}(r,\theta,hD_r,hD_\theta;t,h)
\psi_0(r,\theta)g(r,\theta)^{\frac{1}{4}}
\end{align*}
as operators on ${\mathcal H}$.
Thus
\begin{align*}
&{}\frac{d}{dt}(B_0(t,h)+B_1(t,h))
+i[H,B_0(t,h)+B_1(t,h)]\\
&{}\le g(r,\theta)^{-\frac{1}{4}}
\psi_1(r,\theta)
r_1^{\w}(r,\theta,hD_r,hD_\theta;t,h)\psi_1(r,\theta)
g(r,\theta)^{\frac{1}{4}}.
\end{align*}

We repeat this procedure for $j=2,3,\dots$.
By induction we finally obtain the constants $C_j>0$ 
in the definition of $B_j$ and 
\begin{gather*}
r_j\in S(h^{\mu (j-1)},
\langle h^{-1}t\rangle^{-2}dr^2+d\theta^2
+d\rho^2+d\omega^2),\\
\supp r_j\subset \supp b_j \mod
S(h^\infty)
\end{gather*}
such that 
\begin{align*}
&{}\frac{d}{dt}\sum_{j=0}^kB_j(t,h)
+i\Bigl[H,\sum_{j=0}^kB_j(t,h)\Bigr]\\
&{}\le g(r,\theta)^{-\frac{1}{4}}\psi_k(r,\theta)
r_k^{\w}(r,\theta,hD_r,hD_\theta;t,h)
\psi_k(r,\theta)g(r,\theta)^{\frac{1}{4}},\quad
k=1,2,\dots.
\end{align*}

Let $\psi\in C^\infty(M_\infty)$
be such that it does not depend on $r$ near $\partial M$,
equals $1$ on $\cup_{j=0}^\infty\supp \psi_j$,
and is supported in its sufficiently small neighborhood.
\begin{lemma}\label{14.6.39}
There exists an operator
\begin{align*}
B(t,h)&{}= g(r,\theta)^{-\frac{1}{4}}
\psi(r,\theta)
b^{\w}(r,\theta,hD_r,hD_\theta;t,h)
\psi(r,\theta)
g(r,\theta)^{\frac{1}{4}},\\
b&{}\in S(1,\langle h^{-1}t\rangle^{-2}dr^2
+d\theta^2+d\rho^2+d\omega^2)
\end{align*}
with the following properties:
\begin{enumerate}
\item 
\begin{align*}
B(0,h)=B_0(0,h).
\end{align*}
\item \label{14.6.40}
For any $\varepsilon>0$ we can choose $b$
such that $\supp
b(\cdot,\cdot,\cdot,\cdot;-(t_0+hT_0),h)$
is contained in
\begin{align*}
\{(r,\theta,\rho,\omega)
\in T^*M_\infty;\ 
|r-2h^{-1}\rho^-t_0|<\varepsilon
h^{-1}t_0,\ 
|\theta-\theta_-|<\varepsilon,\\
|\rho-h^{-1}\rho^-|<\varepsilon h^{-1},\ 
|\omega-h^{-1}\omega^-|<\varepsilon h^{-1}\}
\end{align*}
modulo $S(h^\infty)$.
\item The Heisenberg derivative of $B(t,h)$ satisfies
\begin{align*}
\delta B(t,h)
=\frac{d}{dt}B(t,h)+i[H,B(t,h)]
\le
 g(r,\theta)^{-\frac{1}{4}}
\psi(r,\theta)
R(t,h)
\psi(r,\theta)
g(r,\theta)^{\frac{1}{4}}
\end{align*}
on ${\mathcal H}$,
where, through the chart diffeomorphism,
$R(t,h)$ is a bounded operator 
on $L^2(\R^n)$ with
$\sup_{-t_0\le t\le 0}\lVert R(t,h)\rVert=O(h^\infty)$.
\end{enumerate}
\end{lemma}
\textit{Proof of Theorem \ref{13.5.20}.}
We have
\begin{align*}
&{}(e^{-i(t_0+hT_0)H}u_0,B(0,h)e^{-i(t_0+hT_0)H}u_0)_{{\mathcal H}}\\
&{}= ( u_{0},B(-(t_0+hT_0),h)u_{0})_{{\mathcal H}}
\\
&{}\phantom{{}={}}{}
+\int_{-(t_0+hT_0)}^0( e^{-i(t+t_0+hT_0)H}u_0,
\delta B(t,h)e^{-i(t+t_0+hT_0)H}u_0)_{{\mathcal H}}dt\\
&{}
\le O(h^\infty).
\end{align*}
Thus the theorem is proved.
\hfill$\square$


\section{Proof of Main Theorem}
Let $u_0\in{\mathcal H}$ and $t_0>0$.
We suppose $(x_0,\xi^0)\in T^*M\setminus {\mathcal T}_-$, 
and let $(r_-,\theta_-,\rho^-,\omega^-)$ be the 
scattering data at $t=-\infty$.
As in the previous section, 
let $T_0\ll 0$ be large enough that
we can choose a coordinate
neighborhood
$(1,\infty)\times U\subset M_\infty$,
$U\subset \partial M$ such that 
\[\pi\circ\exp tH_p (x_0,\xi^0)
\in (2,\infty)\times U,\quad t\le T_0,\]
and, moreover, 
\begin{align*}
\pi\circ S_{t'}\circ \exp tH_p (x_0,\xi^0)
\in (2,\infty)\times U,\quad t\le T_0,\ t'\le 0.
\end{align*}
By Lemma \ref{7.8.20.9.15} and its analogue
on $M_{\free}$ 
it suffices to show that
\begin{align*}
&{}\exp T_0H_p (x_0,\xi^0)\in
\fs(e^{-i(t_0+hT_0)H}u_0)\\
&{}\iff (r_-+2T_0\rho^-,\theta_-,
\rho^-,\omega^-)
\in \fs(e^{-i(t_0+hT_0)H_0}J^*u_0).
\end{align*}
Let $a\in C^\infty_0(T^*((2,\infty)\times U))$ be a symbol supported 
near $\exp T_0H_p (x_0,\xi^0)$, and 
$\varphi\in C^\infty_0((2,\infty)\times U)$ be a function such that it is 
equal to $1$ on $\cup_{t'\le 0}\pi\circ S_{t'}(\supp a)$ and is 
supported in its small neighborhood. 
We define the operator $A(h)$ by
\begin{align*}
A(h)=\kappa^*
(\kappa_*g)^{-\frac{1}{4}}
(\kappa_*\varphi)
(\kappa_*a)^{\w}(r,\theta,hD_r,hD_\theta)
(\kappa_*\varphi)
(\kappa_*g)^{\frac{1}{4}}
\kappa_*
\quad \mbox{ on } {\mathcal H},
\end{align*}
where $\kappa\colon (1,\infty)\times U
\to (1,\infty)\times \tilde{U}\subset \R^n$ is 
a chart diffeomorphism, and will be
omitted, as in the previous section.
Then, putting $t_0'=t_0+hT_0$ for simplicity, we have
\begin{align}
\begin{split}
\lVert A(h)e^{-it'_0H}u_0\rVert_{{\mathcal H}}^2
={}&( e^{-it'_0H_0}J^*u_0,B(t'_0,h)e^{-it'_0H_0}J^*u_0)_{{\free}}\\
&{}+( e^{-it'_0H}u_0,A(h)^2e^{-it'_0H}(1-JJ^*)u_0)_{{\mathcal H}}\\
&{}+( A(h)^2e^{-it'_0H}(1-JJ^*)u_0,e^{-it'_0H}JJ^*u_0)_{{\mathcal H}},
\end{split}\label{9.10.22.30}
\end{align}
where we set  
\begin{align*}
B(t',h)=e^{-it'H_0}J^*e^{it'H}
A(h)^2e^{-it'H}Je^{it'H_0}
\quad \mbox{ on }{\mathcal H}_{\free}
\end{align*}
for  $t'\ge 0$. We first show that the last two terms in the right-hand side
of (\ref{9.10.22.30}) can be ignored:
\begin{lemma}\label{9.10.22.35}
For any $T_1\ge 0$
\begin{align*}
\lVert A(h)^2e^{-it'H}(1-JJ^*)\rVert_{{\mathcal H}\to {\mathcal H}}
=O(h^\infty)\quad \mbox{ uniformly in }t'\in [0,T_1].
\end{align*}
\end{lemma}
\textit{Proof.}
Let $\tilde{A}(t',h)$ denote the operator $B(-t',h)$
constructed
in the previous section
with $0\le t'\le T_1$. Then 
\begin{align*}
&{}(e^{-it'H}(1-JJ^*)v,\tilde{A}(0,h)e^{-it'H}(1-JJ^*)v)_{{\mathcal H}}\\
&{}\le 
((1-JJ^*)v,\tilde{A}(t',h)(1-JJ^*)v)_{{\mathcal H}}
+T_1\sup_{0\le t'\le T_1}\lVert R(-t',h)\rVert
\lVert (1-JJ^*)v\rVert_{{\mathcal H}}^2.
\end{align*}
Noting the supporting properties of operators
$\tilde{A}(t',h)$ and $1-JJ^*$,
we conclude that the right-hand side is 
$O(h^\infty)$ uniformly in $t'\in [0,T_1]$.
\hfill$\square$

\smallskip
\noindent
In particular, the second and the third terms in 
the right-hand side of (\ref{9.10.22.30}) 
are $O(h^\infty)$.
Thus it suffices to show that
$B(t'_0,h)$ is an $h$-pseudodifferential operator 
with principal symbol $(a\circ S^{-1}_{-t'_0/h})^2$,
since then the principal symbol
approaches $(a\circ S^{-1}_{-\infty})^2$ as $h>0$
tends to $0$.
\begin{lemma}
Let $\tilde{H}$ be the self-adjoint differential 
operator on ${\mathcal H}_{\free}$ defined by
\begin{align*}
\tilde{H}=J^*HJ+(1-J^*J)H_0(1-J^*J).
\end{align*}
Then for any $T_1\ge 0$ we have
\begin{align*}
B(t',h)=e^{-it'H_0}e^{it'\tilde{H}}
(J^*
A(h)^2J)e^{-it'\tilde{H}}e^{it'H_0}
+O(h^\infty)\quad \mbox{ in }
{\mathcal L}({\mathcal H}_{\free})
\end{align*}
uniformly in $t'\in [0,T_1]$.
\end{lemma}
\textit{Proof.}
Similarly to Lemma \ref{9.10.22.35}
we learn by duality that
\begin{align*}
\lVert (1-JJ^*)e^{it'H}A(h)\rVert_{{\mathcal H}\to {\mathcal H}}
=O(h^\infty)
\end{align*}
uniformly in $t'\in [0,T_1]$.
Since
\begin{align*}
(J^*H-\tilde{H}J^*)e^{it'H}A(h)
=(J^*H+(1-J^*J)H_0J^*)(1-JJ^*)e^{it'H}A(h),
\end{align*}
we obtain
\begin{align*}
\lVert e^{-it'\tilde{H}}(J^*H-\tilde{H}J^*)
e^{it'H}A(h)\rVert_{{\mathcal H}\to {\mathcal H}}
=O(h^\infty).
\end{align*}
Hence, by integrating it with respect to $t'$, we have
\begin{align*}
\lVert  J^*e^{it'H}A(h)-e^{it'\tilde{H}}J^*A(h)
\rVert_{{\mathcal H}\to {\mathcal H}}
=O(h^\infty).
\end{align*}
\hfill $\square$

\smallskip
\noindent
Thus we have reduced the proof to the study of the behavior of 
\begin{align*}
\tilde{B}(t',h)=e^{-it'H_0}e^{it'\tilde{H}}(J^*A(h)^2J)
e^{-it'\tilde{H}}e^{it'H_0}.
\end{align*}
We note
\begin{align*}
\frac{d}{dt'}e^{-it'H_0}e^{it'\tilde{H}}
&{}=-i(e^{-it'H_0}(H_0-\tilde{H})e^{it'H_0})
e^{-it'H_0}e^{it'\tilde{H}}\\
&{}=-iL(t')e^{-it'H_0}e^{it'\tilde{H}}.
\end{align*}
If we write $H_0-\tilde{H}=w^{\w}(r,\theta,D_r,D_\theta)$,
then
\begin{align*}
L(t')=w^{\w}(r-2t'D_r,\theta,D_r,D_\theta)
\end{align*}
with no error terms by virtue of the Weyl calculus.
Now $\tilde{B}(t',h)$ satisfies the Heisenberg equation:
\begin{gather*}
\frac{d}{dt}\tilde{B}(t',h)+i[L(t'),\tilde{B}(t',h)]=0,\\
\tilde{B}(0,h)=[h(\theta)^{-\frac{1}{4}}
\varphi(r,\theta)a^{\w}(r,\theta,hD_r,hD_\theta)
\varphi(r,\theta)
h(\theta)^{\frac{1}{4}}]^2.
\end{gather*}
Let us construct an asymptotic solution to this equation.
We put
\begin{align*}
b_0(r,\theta,\rho,\omega;t',h)
=(a\circ S_{-t'/h}^{-1}(r,\theta,\rho,\omega))^2
\end{align*}
for $t'\in [0,T_1]$, $h\in (0,1]$.
Note
\begin{gather*}
b_0(\cdot,\cdot,\cdot,\cdot;t',h)
\in C_0^\infty((2,\infty)\times U),\quad b_0\in S(1),\\
\varphi=1\quad \mbox{ on }\bigcup
\{\pi (\supp b_0(\cdot,\cdot,\cdot,\cdot;t',h));\ 
t'\in [0,T_1],\  h\in (0,1]\}.
\end{gather*}
Since we are considering operators on $M_{\free}$,
we may use chart neighborhoods of the form
$\R\times U$, $U\subset M_\infty$
in what follows.
Choose a cutoff function $\psi_0\in C^\infty(\R\times U)$ 
which does not depend on $r$,
equals $1$ on $\supp \varphi$ and
is supported in its sufficiently small neighborhood,
and define
\begin{align*}
B_0(t',h)=h(\theta)^{-\frac{1}{4}}
\psi_0(\theta)
b_0^{\w}(r,\theta,hD_r,hD_\theta;t',h)
\psi_0(\theta)h(\theta)^{\frac{1}{4}}.
\end{align*}
\begin{lemma}\label{9.20.9.50}
There is $r_0\in S(\langle h^{-1}t'\rangle^{-1-\mu})$ which is 
supported in $\supp b_0$ modulo $S(h^\infty)$ such that
\begin{align*}
&{}\frac{d}{dt'}B_0(t',h)+i[L(t'),B_0(t',h)]\\
&{}=h(\theta)^{-\frac{1}{4}}
\psi_1(\theta)
r_0^{\w}(r,\theta,hD_r,hD_\theta ;t',h)
\psi_1(\theta)h(\theta)^{\frac{1}{4}}.
\end{align*}
Here $\psi_1\in C^\infty(\R\times U)$ is any
function that does not depend on $r$ and
equals $1$ on $\supp \psi_0$.
\end{lemma}
\textit{Proof.}\ 
Differentiating 
\begin{align*}
b_0(S_{-t'/h}(r,\theta,\rho,\omega);t',h)
=a(r,\theta,\rho,\omega)^2
\end{align*}
with respect to $t'$, we obtain
\begin{align*}
\frac{\partial }{\partial t'}b_0+h\{k,b_0\}
=0,\quad 
b_0(r,\theta,\rho,\omega;0,h)
=a(r,\theta,\rho,\omega)^2,
\end{align*}
where $k$
is the principal part of
$w(r-2h^{-1}t'\rho, \theta,h^{-1}\rho,h^{-1}\omega)$:
\begin{align*}
k(r,\theta,\rho,\omega;t',h)=j(r-2h^{-1}t'\rho)^2[h^{-2}\rho^2-
p(r-2h^{-1}t'\rho,\theta,h^{-1}\rho,h^{-1}\omega)].
\end{align*}
Since $r,\rho,\omega=O(1)$ and $\rho<c<0$
for
$(r,\theta,\rho,\omega,t',h)\in \supp b_0$, we have 
\begin{align}
|\partial_r^j\partial_\theta^\alpha
\partial_\rho^k\partial_\omega^\beta
w(r-2h^{-1}t'\rho, \theta,h^{-1}\rho,h^{-1}\omega)|
\le Ch^{-2}\langle h^{-1}t'\rangle^{-1-\mu-j}
\quad \mbox{ on }\supp b_0,
\label{10.10.8.15}
\end{align}
and thus, if we apply the asymptotic expansion formula,
\begin{align*}
&{}\frac{d}{dt'}B_0(t',h)+i[L(t'),B_0(t',h)]\\
&{}=h(\theta)^{-\frac{1}{4}}
\psi_1(\theta)
\Bigl[\sum_{j=0}^{N-1}h^{j} r_{0,j}+R_N\Bigr]
^{\w}(r,\theta,hD_r,hD_\theta ;t',h)
\psi_1(\theta)h(\theta)^{\frac{1}{4}}
\end{align*}
with
\begin{align*}
r_{0,j}\in S(\langle h^{-1}t'\rangle^{-1-\mu}).
\end{align*}
Therefore
it suffices to confirm that the remainder term
$R_N$ has an appropriate decaying property in $h$, 
since, then, 
the required $r_0$ will be obtained
by taking the asymptotic sum of $r_{0,j}$.
However, since the estimate (\ref{10.10.8.15})
does hold on the whole $\R^{2n}$,
we can not estimate it directly.
Let us rewrite
\begin{align*}
i[L(t'), B_0(t',h)]
&{}=ie^{-itH_0}[H_0-\tilde{H},e^{itH_0}B_0(t',h)e^{-itH_0}]e^{itH_0},\\
e^{itH_0}B_0(t',h)e^{-itH_0}&{}=h(\theta)^{-\frac{1}{4}}
\psi_0(\theta)
b_0^{\w}(r+2t'D_r,\theta,hD_r,hD_\theta;t',h)
\psi_0(\theta)h(\theta)^{\frac{1}{4}},
\end{align*}
and compute the remainder term of
$h(\theta)^{\frac{1}{4}}[H_0-\tilde{H},
e^{itH_0}B_0(t',h)e^{-itH_0}]h(\theta)^{-\frac{1}{4}}$,
instead.
The $N$-th remainder term 
is given by
\begin{align}
\begin{split}
{\mathcal R}_N[&w(r, \theta,h^{-1}\rho,h^{-1}\omega)
b_0(r_1+2h^{-1}t'\rho_1, \theta_1,\rho_1,\omega_1;t',h)\\
&{}-b_0(r+2h^{-1}t'\rho, \theta,\rho,\omega;t',h)
w(r_1, \theta_1,h^{-1}\rho_1,h^{-1}\omega_1)
]\bigg|_{(r_1,\theta_1,\rho_1,\omega_1)
=(r,\theta,\rho,\omega)}
\end{split}\label{9.26.21.45}
\end{align}
modulo $S(h^\infty)$ (see, e.g., \cite{Hormander}).
Here
\begin{align*}
{\mathcal R}_N={}&
\int_0^1\frac{(1-\tau)^{N-1}}{(N-1)!}
e^{\frac{ih}{2}\tau(D_rD_{\rho_1}+D_{\theta}D_{\omega_1}
-D_\rho D_{r_1}-D_{\omega}D_{\theta_1})}d\tau \\
{}&\cdot\Bigl\{\frac{ih}{2}(
\partial_r\partial_{\rho_1}+\partial_{\theta}\partial_{\omega_1}
-\partial_\rho \partial_{r_1}-\partial_{\omega}\partial_{\theta_1}
)\Bigr\}^N,
\end{align*}
which is defined as a
Fourier multiplier by
\begin{align*}
\int_0^1\frac{(1-\tau)^{N-1}}{(N-1)!}
e^{\frac{ih}{2}\tau(\hat{r}\hat{\rho}_1
+\hat{\theta}\hat{\omega}_1
-\hat{\rho}\hat{r}_1-\hat{\omega}\hat{\theta}_1)}d\tau 
\Bigl\{\frac{ih}{2}(
\hat{r}\hat{\rho}_1
+\hat{\theta}\hat{\omega}_1
-\hat{\rho}\hat{r}_1
-\hat{\omega}\hat{\theta}_1)\Bigr\}^N.
\end{align*}
Write (\ref{9.26.21.45}) in 
the integral form with respect to
$(r,\theta,\rho,\omega,r_1,\theta_1,\rho_1,\omega_1)$
and their conjugate variables
$(\hat{r},\hat{\theta},\hat{\rho},\hat{\omega},
\hat{r}_1,\hat{\theta}_1,\hat{\rho}_1,\hat{\omega}_1)$,
and integrate it by parts to make it integrable.
By the supporting property of $b_0$
we obtain, for example,
\begin{gather*}
\partial_r^j\partial_\theta^\alpha
\partial_\rho^k\partial_\omega^\beta
b_0(r+2h^{-1}t'\rho, \theta,\rho,\omega;t',h)
\in S(\langle r\rangle^{k},
dr^2+d\theta^2+h^{-2}t'{}^2d\rho^2+d\omega^2),\\
w(r, \theta,h^{-1}\rho,h^{-1}\omega)
\in S(h^{-2}\langle \rho;\omega\rangle^2,
\langle r\rangle^{-2}dr^2+d\theta^2
+\langle \rho\rangle^{-2}d\rho^2+\langle \omega\rangle^{-2}d\omega^2),
\end{gather*}
and we learn that (\ref{9.26.21.45}) belongs to
\begin{align*}
 S(h^{N-M},
dr^2+d\theta^2
+h^{-2}d\rho^2+d\omega^2),
\end{align*}
where $M$ is chosen independent of $N$.
Then it is easy to see that the asymptotic expansion is verified.
\hfill $\square$

\medskip
Now we solve the second transport equation:
\begin{align*}
\frac{\partial}{\partial t'}b_1+h\{k,b_1\}
=-r_0,\quad
b_1(r,\theta,\rho,\omega;0,h)=0.
\end{align*}
The solution is given by
\begin{align*}
b_1(r,\theta,\rho,\omega;t',h)
=\int^{t'}_0r_0(S_{-s/h}\circ
S_{-t'/h}^{-1}(r,\theta,\rho,\omega); s,h)\,ds,
\end{align*}
and it
satisfies by Lemma \ref{9.20.9.50}
\begin{align*}
b_1\in S(h), \quad
\supp b_1\subset \supp b_0 \mod S(h^\infty).
\end{align*}
Similarly to Lemma \ref{9.20.9.50} 
there is 
$r_1\in S(h\langle h^{-1}t'\rangle^{-1-\mu})$ 
such that  it has 
support in $\supp b_0$ modulo $S(h^\infty)$ and  that
\begin{align*}
&{}\frac{d}{dt'}(B_0(t',h)+B_1(t',h))
+i[L(t'),B_0(t',h)+B_1(t',h)]\\
&{}=h(\theta)^{-\frac{1}{4}}
\psi_2(\theta)
r_1^{\w}(r,\theta,hD_r,hD_\theta ;t',h)
\psi_2(\theta)h(\theta)^{\frac{1}{4}},
\end{align*}
where
\begin{align*}
B_1(t',h)=h(\theta)^{-\frac{1}{4}}
\psi_1(\theta)
b_1^{\w}(r,\theta,hD_r,hD_\theta;t',h)
\psi_1(\theta)h(\theta)^{\frac{1}{4}}
\end{align*}
and $\psi_2\in C^\infty(\R\times U)$ is any
function that does not depend on $r$ and
equals $1$ on $\supp \psi_1$.

Repeating this procedure, we can inductively
construct the
symbols $b_j$ and $r_j$.
Let $\psi\in C^\infty(\R\times U)$ be 
independent of $r$ and equal to $1$ on
$\cup_{j=0}^\infty\supp \psi_j$,
and set 
\begin{align*}
\Tilde{\Tilde{B}}(t',h)=h(\theta)^{-\frac{1}{4}}
\psi(\theta)
b^{\w}(r,\theta,hD_r,hD_\theta;t',h)
\psi(\theta)h(\theta)^{\frac{1}{4}},\quad 
b\sim\sum_{j=0}^\infty b_j,
\end{align*}
then we have 
\begin{gather*}
\frac{d}{d t'}\Tilde{\Tilde{B}}(t',h)
+i[L(t'),\Tilde{\Tilde{B}}(t',h)]=O(h^\infty),\\
\Tilde{\Tilde{B}}(0,h)=[h(\theta)^{-\frac{1}{4}}
\varphi(r,\theta)a^{\w}(r,\theta,hD_r,hD_\theta)
\varphi(r,\theta)
h(\theta)^{\frac{1}{4}}]^2.
\end{gather*}
This implies
\begin{align*}
B(t',h)=\tilde{B}(t',h)+O(h^\infty)
=\Tilde{\Tilde{B}}(t',h)+O(h^\infty).
\end{align*}
Substituting this to (\ref{9.10.22.30}), 
we obtain
\begin{align*}
\lVert A(h)e^{-it'_0H}u_0\rVert_{{\mathcal H}}^2
=( e^{-it'_0H_0}J^*u_0,
\Tilde{\Tilde{B}}(t'_0,h)e^{-it'_0H_0}J^*u_0)_{\free}
+O(h^\infty).
\end{align*}
$\Tilde{\Tilde{B}}(t'_0, h)$ is
an $h$-pseudodifferential operator with principal symbol
$(a\circ S_{-t_0'/h}^{-1})^2$.

Assume $(r_-,\theta_-,\rho^-,\omega^-)
\notin \wf(e^{-it_0H_0}J^*u_0)$,
or equivalently
\begin{align*}
(r_-+2T_0\rho^-,\theta_-,
\rho^-,\omega^-)
\notin \fs(e^{-i(t_0+hT_0)H_0}J^*u_0).
\end{align*}
If the symbol $a$ is supported sufficiently near
$\exp T_0H_p (x_0,\xi^0)$, then the full symbol of 
$\Tilde{\Tilde{B}}(t_0+hT_0, h)$ is
supported near $(r_-+2T_0\rho^-,\theta_-,
\rho^-,\omega^-)$
for small $h>0$,
so that
\begin{align*}
( e^{-i(t_0+hT_0)H_0}J^*u_0,
\Tilde{\Tilde{B}}(t_0+hT_0,h)e^{-i(t_0+hT_0)H_0}J^*u_0)_{\free}=
O(h^\infty).
\end{align*}
Thus we obtain
\begin{align*}
\lVert A(h)e^{-i(t_0+hT_0)H}u_0\rVert_{{\mathcal H}}=O(h^\infty),
\end{align*}
which implies
$\exp T_0H_p (x_0,\xi^0)\notin
\fs(e^{-i(t_0+hT_0)H}u_0)$.
Hence it follows that
$(x_0,\xi^0)\notin \wf(e^{-it_0H}u_0)$.

Conversely, assume $(x_0,\xi^0)\notin \wf(e^{-it_0H}u_0)$.
We choose a symbol $a\ge 0$
supported in a small neighborhood of
$\exp T_0H_p (x_0,\xi^0)$
and construct the operator
$\Tilde{\Tilde{B}}(t_0+hT_0,h)$ as above accordingly.
Then by the assumption we have
\begin{align*}
( e^{-i(t_0+hT_0)H_0}J^*u_0,
\Tilde{\Tilde{B}}(t_0+hT_0,h)e^{-i(t_0+hT_0)H_0}J^*u_0)_{\free}=
O(h^\infty).
\end{align*}
The principal symbol of 
$\Tilde{\Tilde{B}}(t_0+hT_0, h)$ is
positive near $(r_-+2T_0\rho^-,\theta_-,\rho^-,\omega^-)$
for small $h>0$, and thus it follows that
\begin{align*}
(r_-+2T_0\rho^-,\theta_-,
\rho^-,\omega^-)
\notin \fs(e^{-i(t_0+hT_0)H_0}J^*u_0),
\end{align*}
or equivalently $(r_-,\theta_-,\rho^-,\omega^-)
\notin \wf(e^{-it_0H_0}J^*u_0)$.
The proof of Theorem \ref{7.8.20.1.45} is complete.

\end{document}